\input amstex
\input amsppt.sty
\magnification=\magstep1
\hsize=30truecc
\vsize=22.2truecm
\baselineskip=16truept
\TagsOnRight
\nologo
\pageno=1
\topmatter
\def\N{\Bbb N}
\def\Z{\Bbb Z}
\def\Q{\Bbb Q}

\def\C{\Bbb C}
\def\l{\left}
\def\r{\right}
\def\b{\bigg}
\def\bg{\bigg}
\def\({\b(}
\def\[{\b[}
\def\){\b)}
\def\]{\b]}

\def\t{\text}
\def\f{\frac}
\def\mo{\roman{mod}}
\def\ord{\roman{ord}}

\def\sm{\setminus}

\def\bi{\binom}
\def\eq{\equiv}

\def\ls{\leqslant}
\def\gs{\geqslant}
\def\al{\alpha}
\def\ve{\varepsilon}

\def\Proof{\noindent{\it Proof}}
\def\Remark{\medskip\noindent{\it Remark}}

\hbox {Acta Arith. {\bf 127}(2007), no.\,4, 337--363.}
\bigskip
\title On Fleck quotients\endtitle
\author  Zhi-Wei Sun$^1$ (Nanjing) and Daqing Wan$^2$ (Irvine, CA) \endauthor
\leftheadtext{Zhi-Wei Sun and Daqing Wan}
\affil $^1$Department of Mathematics, Nanjing University
\\ Nanjing 210093, People's Republic of China
\\zwsun\@nju.edu.cn
\\ {\tt http://math.nju.edu.cn/$^\sim$zwsun}
\medskip
$^2$Department of Mathematics, University of California\\Irvine,
CA 92697-3875, USA
\\dwan\@math.uci.edu
\\ {\tt http://www.math.uci.edu/$^{\sim}$dwan}
\medskip
\endaffil
\thanks 2000 {\it Mathematics Subject Classification}.\,Primary 11B65;
Secondary 05A10, 11A07, 11B68, 11B73, 11L05, 11S99.\newline\indent
The first author is supported by the National Science Fund for
Distinguished Young Scholars (No. 10425103) in China. The
second author is partially supported by NSF.
\endthanks
\abstract Let $p$ be a prime, and let $n\gs1$ and $r$ be integers.
In this paper we study Fleck's quotient
$$F_p(n,r)=(-p)^{-\lfloor(n-1)/(p-1)\rfloor}
\sum_{k\eq r\,(\mo\ p)}\bi nk(-1)^k\in\Z.$$
We determine $F_p(n,r)$ mod $p$ completely by certain number-theoretic
and combinatorial methods; consequently, if $2\ls n\ls p$ then
$$\sum_{k=1}^n(-1)^{pk-1}\bi{pn-1}{pk-1}\eq(n-1)!B_{p-n}p^n\ (\mo\ p^{n+1}),$$
where $B_0,B_1,\ldots$ are Bernoulli numbers.
We also establish the Kummer-type congruence
$F_p(n+p^a(p-1),r)\eq F_p(n,r)\ (\mo\ p^a)$ for $a=1,2,3,\ldots$,
and reveal some connections between Fleck's quotients and class numbers
of the quadratic fields $\Q(\sqrt{\pm p})$
and the $p$-th cyclotomic field $\Q(\zeta_p)$. In addition,
generalized Fleck quotients are also studied in this paper.
\endabstract
\endtopmatter

\document

\heading{1. Introduction and main results}\endheading

Let $m\in\Z^+=\{1,2,\ldots\}$, $n\in\N=\{0,1,\ldots\}$ and $r\in\Z$,
and define
$$C_m(n,r)=\sum_{k\eq r\, (\mo\ m)}\bi nk(-1)^k.\tag1.0$$
This sum has been studied by various authors and many applications have been found
(cf. [S02] and its references).
The following well-known observation is fundamental:
$$mC_m(n,r)=\sum_{k=0}^n\bi nk(-1)^k\sum_{\gamma^m=1}\gamma^{k-r}
=\sum_{\gamma^m=1}\gamma^{-r}(1-\gamma)^n.$$
Note that
$$C_m(n+1,r)=C_m(n,r)-C_m(n,r-1)$$
since $x^{-r}(1-x)^{n+1}=x^{-r}(1-x)^n-x^{-r+1}(1-x)^n$.

Let $p$ be a prime, and let $n\in\N$ and $r\in\Z$.
In 1913 A. Fleck (cf. [D, p. 274]) showed that
$$\ord_p(C_p(n,r))\gs\l\lfloor\f{n-1}{p-1}\r\rfloor,$$
where $\ord_p(\al)$ denotes the $p$-adic order
of a $p$-adic number $\al$, and $\lfloor\cdot\rfloor$
is the well-known floor function.
Fleck's result is fundamental
in the recent investigation of the $\psi$-operator related
to Fontaine's theory, Iwasawa's theory, and $p$-adic Langlands
correspondence (cf. [Co], [SW] and [W]); it also plays
an indispensable role in Davis and Sun's study of
homotopy exponents of special unitary groups
(cf. [DS] and [SD]).
In this paper we are interested in the {\it Fleck quotient}
$$F_p(n,r):=(-p)^{-\lfloor(n-1)/(p-1)\rfloor}
C_p(n,r)+[\![n=0]\!].\tag1.1$$
(Throughout this paper, for an assertion $A$ we let $[\![A]\!]$ take $1$ or $0$
according as $A$ holds or not.)

For $a\in\Z$ and $m\in\Z^+$, we use $\{a\}_m$
to denote the least nonnegative residue of $a$ mod $m$
(thus $\{a\}_m/m$ is the fractional part $\{a/m\}$
of $a/m$).
For a prime $p$ and an integer $a$, we define $q_p(a)=(a^{p-1}-1)/p$
which is an integer if $a\not\eq0\ (\mo\ p)$.

By a number-theoretic approach related to Gauss sums,
we establish the following explicit result.
\proclaim{Theorem 1.1} Let $p$ be a prime, and let $n\in\N$ and $r\in\Z$.
Set $n_0=\{n\}_p$ and $n_1=\{n_0-n\}_{p-1}=\{-\lfloor n/p\rfloor\}_{p-1}$.
If $n_0\ls n_1$, then
$$F_p(n,r)\eq\f{(-1)^{n_1}}{n_1!}
\sum_{k=0}^{n_0}\bi{n_0}k(-1)^k(k-r)^{n_1}\ (\mo\ p).\tag1.2$$
If $n_0>n_1=0$, then
$$F_p(n,r)\eq(-1)^{\{r\}_p}\bi{n_0}{\{r\}_p}\ (\mo\ p).\tag1.3$$
If $n_0>n_1>0$, then
$$F_p(n,r)\eq\f{(-1)^{n_1-1}}{(n_1-1)!}
\sum_{k=0}^{n_0}\bi{n_0}k(-1)^{k}(k-r)^{n_1}q_p(k-r)\ (\mo\ p).
\tag1.4$$
\endproclaim

\proclaim{Corollary 1.1} Let $p$ be a prime and let $n\in\N$ and $r\in\Z$. Then
$$F_p(pn,r)\eq\f{r^{n^*}}{n^*!}\ (\mo\ p)\tag1.5$$
where $n^*=\{-n\}_{p-1}$. Consequently,
$$F_p\l(p\f{p-1}2,r\r)\eq
\cases(-1)^{(h(-p)+1)/2}(\f rp)\ (\mo\ p)&\t{if}\ p\not=3\ \&\ 4\mid p+1,
\\(-1)^{(h(p)-1)/2}(\f rp)\f v2\ (\mo\ p)&\t{if}\ 4\mid p-1,
\endcases\tag1.6$$
where $(\f{\cdot}p)$ is the Legendre symbol, and $h(-p)$ and $h(p)$
are the class numbers of the quadratic fields $\Q(\sqrt{-p})$ and $\Q(\sqrt{p})$
respectively, and for $p\eq1\ (\mo\ 4)$ we write the fundamental unit
of $\Q(\sqrt p)$ in the form $(v+u\sqrt p)/2$ with $u,v\in\Z$ and $u\eq v\ (\mo\ 2)$.
\endproclaim
\Proof. Note that $\{pn\}_p=0$.
By Theorem 1.1,
$$F_p(pn,r)\eq\f{(-1)^{n^*}}{n^*!}\sum_{k=0}^0\bi0k(-1)^k(k-r)^{n^*}
=\f{r^{n^*}}{n^*!}\ (\mo\ p).$$

When $p\not=2$ and $n=(p-1)/2$, we have $n^*=(p-1)/2$ and hence
$$\align F_p\l(p\f{p-1}2,r\r)
\eq&r^{(p-1)/2}(-1)^{(p-1)/2}\f{((p-1)/2)!}{\prod_{k=1}^{(p-1)/2}k(p-k)}
\\\eq&\l(\f rp\r)(-1)^{(p-1)/2}\f{((p-1)/2)!}{(p-1)!} \  (\t{by Euler's criterion})
\\\eq&(-1)^{(p+1)/2}\l(\f rp\r)\f{p-1}2!\ (\mo\ p)\, (\t{by Wilson's theorem}).
\endalign$$
If $p>3$ and $p\eq 3\ (\mo\ 4)$, then
$$\f{p-1}2!\eq(-1)^{(h(-p)+1)/2}\ (\mo\ p)$$
by a result of  L. J. Mordell [M].
When $p\eq 1\ (\mo\ 4)$ and $\ve_p=(v+u\sqrt p)/2>1$ is the fundamental unit
of $\Q(\sqrt p)$ with $u,v\in\Z$ and $u\eq v\ (\mo\ 2)$, by S. Chowla [C] we have
$$\f{p-1}2!\eq(-1)^{(h(p)+1)/2}\,\f v2\ (\mo\ p).$$
Combining the above we immediately obtain (1.6). \qed

\Remark. Let $n$ be a positive integer and $p>2n+1$ be a prime.
By the first part of Corollary 1.1
in the case $r=0$, we have
$$\bi{2pn}{pn}(-1)^n+2\sum_{k=0}^{n-1}\bi{2pn}{pk}(-1)^k
=\sum_{k=0}^{2n}\bi{2pn}{pk}(-1)^{pk}
\eq0\ (\mo\ p^{2n+1})$$
and hence
$$\bi{2pn-1}{pn-1}=\f12\bi{2pn}{pn}
\eq\sum_{k=0}^{n-1}(-1)^{n-1-k}\bi{2pn}{pk}\ (\mo\ p^{2n+1}).\tag1.7$$
When $n=1$ and $p>3$, this
gives the Wolstenholme congruence
$$\f12\bi{2p}p=\bi{2p-1}{p-1}\eq1\ (\mo\ p^3).$$
When $n=2$ and $p>5$, (1.7)
yields the following new congruence
$$\bi{4p-1}{2p-1}=\f12\bi{4p}{2p}\eq\bi{4p}p-1\ \l(\mo\ p^5\r).$$
\medskip

Our second approach to Fleck quotients is of combinatorial nature.
It involves Stirling numbers of the second kind as well as
higher-order Bernoulli polynomials.

Let $n\in\N$. The Stirling numbers $S(n,k)\ (k\in\N)$
of the second kind are given by
$$x^n=\sum_{k\in\N}S(n,k)(x)_k,$$
where
$$(x)_0=1\ \ \t{and}\ \ (x)_k=x(x-1)\cdots(x-k+1)\ \t{for}\ k=1,2,\ldots.$$
Clearly, $S(n,n)=1$, and $S(n,k)=0$ if $k>n$.
When $n+k>0$,
$S(n,k)$ is actually the number of ways to partition a set of cardinality $n$
into $k$ nonempty subsets. Here is an explicit formula (cf. [LW, p.\,126]) for
Stirling numbers of the second kind:
$$S(n,k)=\f1{k!}\sum_{j=0}^k\bi kj(-1)^{k-j}j^n.$$
As $S(i,k)=0$ for all those $i\in\N$ with $i<k$, we have {\it Euler's identity}
$$\sum_{j=0}^k\bi kj(-1)^{j}P(j)=0,$$
where $P(x)$ is any polynomial with $\deg P<k$ having complex number coefficients.
It is known (cf. [LW, p.\,126]) that
$$\sum_{n=k}^{\infty}S(n,k)\f{x^n}{n!}=\f{(e^x-1)^k}{k!};$$
in other words,
$$(e^x-1)^k=\sum_{n=k}^{\infty}\bar S(n,k)x^n \ \ \t{with}\ \ \bar S(n,k)=\f{k!}{n!}S(n,k).$$

For $m=0,1,\ldots$, the $m$-th order Bernoulli polynomials $B_n^{(m)}(t)\ (n\in\N)$
are defined by
$$\f{x^me^{tx}}{(e^x-1)^m}=\sum_{n=0}^{\infty}B_n^{(m)}(t)\f{x^n}{n!},\tag1.8$$
and those $B_n^{(m)}=B_n^{(m)}(0)$ are called the $m$-th order Bernoulli numbers.
The usual Bernoulli polynomials and numbers are $B_n(t)=B_n^{(1)}(t)$
and $B_n=B_n(0)=B_n^{(1)}$ respectively.
(It is well known that $B_0=1$, $B_1=-1/2$ and $B_{2k+1}=0$ for $k=1,2,\ldots$;
the reader may consult [IR, pp.\,228--248] for the basic properties of Bernoulli numbers.)
For a formal power series
$f(x)=\sum_{n=0}^{\infty}a_nx^n$, we use $[x^n]f(x)$
to denote the coefficient $a_n$ of the monomial $x^n$ in $f(x)$.
Thus
$$\align B_n^{(m)}(t)=&[x^n]n!\l(\f{x}{e^x-1}\r)^me^{tx}
\\=&[x^n]n!\sum_{k=0}^{\infty}B_k^{(m)}\f{x^k}{k!}\sum_{j=0}^{\infty}\f{(tx)^j}{j!}
=\sum_{k=0}^n\bi nkB_k^{(m)}t^{n-k}.
\endalign$$
It is also easy to verify that
$B_n^{(m)}(m-t)=(-1)^nB_n^{(m)}(t)$, and
$$\f{B_n^{(m)}(t)}{n!}
=\sum_{k_0+\cdots+k_{m-1}=n}\f{B_{k_0}(t)}{k_0!}\prod_{0<i<m}\f{B_{k_i}}{k_i!}\quad \t{provided}\ m>0.$$
If $0\ls n<p-1$, then $B_0,\ldots,B_n$ are $p$-adic integers by
the von Staudt-Clausen theorem (cf. [IR, p.\,233])
or the recurrence $\sum_{k=0}^l\bi{l+1}kB_k=0\ (l=1,2,\ldots)$,
therefore $B_n^{(m)}(t)\in\Z_p[t]$ where $\Z_p$ is the ring of $p$-adic integers.

Our discovery of the next theorem was actually motivated by Theorem 1.1.

\proclaim{Theorem 1.2} Let $p$ be a prime, and let $n\in\N$ and $r\in\Z$.
Set $n^*=\{-n\}_{p-1}$.
For any integer $m\eq n\ (\mo\ p)$, if $m\gs0$ then
$(-1)^nF_p(n,r)$ is congruent to
$$\aligned \sum_{k=0}^{n^*}\bar S(n^*-k+m,m)\f{(-r)^k}{k!}
=&\sum_{k=0}^{n^*}\bar S(m+n^*,m+k)\bi{-r}k
\\=&\sum_{k=0}^m\bi mk(-1)^{m-k}\f{(k-r)^{m+n^*}}{(m+n^*)!}
\endaligned\tag1.9$$
modulo $p$;
if $m\ls 0$ then we have
$$F_p(n,r)\eq \f{(-1)^{n^*}}{n^*!}B_{n^*}^{(-m)}(-r)
\eq-(p-1-n^*)!B_{n^*}^{(-m)}(-r)\ (\mo\ p).\tag1.10$$
\endproclaim

The following consequence determines $B_n^{(m)}(a)$ modulo a prime $p$
for $m\in\{1,\ldots,p\}$, $n\in\{0,\ldots,p-2\}$ and $a\in\Z$.

\proclaim{Corollary 1.2}
Let $p$ be a prime and $r\in\Z$.
Let $n_0\in\{0,\ldots,p-1\}$ and $n_1\in\{0,\ldots,p-2\}$.
If $n_0\ls n_1$, then
$$B_{n_1-n_0}^{(p-n_0)}(-r)
\eq\f{1}{(n_1)_{n_0}}\sum_{k=0}^{n_0}\bi{n_0}k(-1)^{n_0-k}(k-r)^{n_1}\ (\mo\ p).
\tag1.11$$
If $n_0>n_1=0$, then
$$B_{p-n_0+n_1-1}^{(p-n_0)}(-r)
\eq\f{(-1)^{\{r\}_p-1}}{n_0!}\bi{n_0}{\{r\}_p}\ (\mo\ p).\tag1.12$$
If $n_0>n_1>0$, then
$$\aligned
B_{p-n_0+n_1-1}^{(p-n_0)}(-r)
\eq&\f{(-1)^{n_1}}{(n_0-n_1)!(n_1-1)!}
\\&\times\sum_{k=0}^{n_0}\bi{n_0}k(-1)^{k}(k-r)^{n_1}q_p(k-r)\ (\mo\ p).
\endaligned\tag1.13$$
\endproclaim

\Proof. Let $n$ be a nonnegative integer with $n\eq n_0-pn_1\ (\mo\ p(p-1))$.
Applying (1.10) with $m=n_0-p$ we obtain
$$F_p(n,r)\eq\f{(-1)^{n^*}}{n^*!}B_{n^*}^{(p-n_0)}(-r)
\eq-(p-1-n^*)!B_{n^*}^{(p-n_0)}(-r)\ (\mo\ p),$$
where $n^*=\{-n\}_{p-1}$.

If $n_0\ls n_1$, then $n^*=n_1-n_0$ and hence
$$B_{n_1-n_0}^{(p-n_0)}(-r)\eq(-1)^{n_1-n_0}(n_1-n_0)!F_p(n,r)\ (\mo\ p),$$
which implies (1.11) with the help of (1.2).

Now we consider the case $n_0>n_1$. Clearly $n^*=n_1-n_0+p-1$ and $p-1-n^*=n_0-n_1$.
Therefore
$$F_p(n,r)\eq-(n_0-n_1)!B_{n_1-n_0+p-1}^{(p-n_0)}(-r)\ (\mo\ p).$$
The case $n_1=0$ of this, together with (1.3), yields (1.12).
When $n_1>0$, combining the last congruence with (1.4) we obtain (1.13). \qed

\proclaim{Corollary 1.3} Let $p$ be a prime and let $n\in\Z^+$.
Then $\ord_p(C_p(n,r))=\lfloor(n-1)/(p-1)\rfloor$
for at least $p-n^*\gs 2$ values of $r\in\{0,\ldots,p-1\}$, where $n^*=\{-n\}_{p-1}$.
\endproclaim
\Proof. For any $r\in\Z$,  $\ord_p(C_p(n,r))=\lfloor(n-1)/(p-1)\rfloor$
if and only if $F_p(n,r)\not\eq0\ (\mo\ p)$.
By Theorem 1.2,
$$F_p(n,r)\eq\f{(-1)^{n^*}}{n^*!}B_{n^*}^{(p-\{n\}_p)}(-r)\ (\mo\ p)
\quad \ \t{for all}\ r=0,\ldots,p-1.$$
Recall that $B_{n^*}^{(p-\{n\}_p)}(x)\in\Z_p[x]$
is monic and of degree $n^*$. Also, a polynomial of degree $n^*$ over the field
$\Z/p\Z$ cannot have more than $n^*$ distinct zeroes in the field (cf. [IR, p.\,39]).
So the congruence equation $F_p(n,r)\eq0\ (\mo\ p)$ has at most $n^*$ solutions
with $r\in\{0,\ldots,p-1\}$. This yields the desired result. \qed

\proclaim{Corollary 1.4} Let $p$ be a prime, and let $n\in\N$ and $n^*=\{-n\}_{p-1}$.
Then
$$(-1)^nF_p(n,0)\eq\bar S(n^*+\{n\}_p,\{n\}_p)
\eq\f{B_{n^*}^{(m)}}{n^*!}\ (\mo\ p),\tag1.14$$
where $m$ is any nonnegative integer with $m+n\eq0\ (\mo\ p)$.
Also,
$$(-1)^nF_p(pn+p-1,r)\eq\f{B_{n^*}(-r)}{n^*!}\eq-(p-1-n^*)!B_{n^*}(r+1)\ (\mo\ p)\tag1.15$$
for all $r\in\Z$, and in particular
$$\bi{2p-1}{p+r}+(-1)^p\bi{2p-1}r\eq(-1)^rp^2B_{p-2}(-r)\ (\mo\ p^3)\tag1.16$$
for every $r=0,\ldots,p-1$.
\endproclaim
\Proof.  Applying Theorem 1.2 with $r=0$ we immediately get (1.14).

As $pn+p-1\eq-1\ (\mo\ p)$ and $n^*=\{-(pn+p-1)\}_{p-1}$,
by the second part of Theorem 1.2 and the identity $(-1)^{n^*}B_{n^*}(x)=B_{n^*}(1-x)$,
whenever $r\in\Z$ we have
$$\align(-1)^{n^*}F_p(pn+p-1,r)\eq&\f{B_{n^*}(-r)}{n^*!}
\eq(-1)^{n^*+1}(p-1-n^*)!B_{n^*}(-r)
\\\eq&-(p-1-n^*)!B_{n^*}(r+1)
\ (\mo\ p)\endalign$$
and hence (1.15) holds.

Now let $r\in\{0,\ldots,p-1\}$. By (1.15) in the case $n=1$,
$$-F_p(2p-1,r)\eq-(p-1-(p-2))!B_{p-2}(r+1)\ (\mo\ p)$$
and hence
$$F_p(2p-1,r)\eq B_{p-2}(1-(-r))=(-1)^{p-2}B_{p-2}(-r)\ (\mo\ p)$$
which is equivalent to (1.16). We are done. \qed

\medskip

Let $p$ be an odd prime, and let $h_p$ and $h_p^+$
denote the class numbers of the cyclotomic field $\Q(\zeta_p)$
and its maximal real subfield $\Q(\zeta_p+\zeta_p^{-1})$ respectively,
where $\zeta_p$ is a primitive $p$-th root of unity in the complex field $\C$.
It is well known that $h_p^-=h_p/h_p^+$ is an integer.
If $p$ divides none of
the numerators of the Bernoulli numbers $B_0, B_2,\ldots,B_{p-3}\in\Z_p$,
then $p$ is said to be a {\it regular} prime.
In 1850 E. Kummer proved that
$$\align p\nmid h_p&\iff p\nmid h_p^-\iff p\ \t{is regular}
\\&\ \Longrightarrow \ x^p+y^p=z^p\ \t{has no integer solution with}\ xyz\not=0.
\endalign$$
Furthermore,
$$h_p^-\eq\prod_{0<n\ls(p-3)/2}\l(-\f{B_{2n}}{4n}\r)\ (\mo\ p)$$
by the proof of Theorem 5.16 in [Wa, p.\,62].

\proclaim{Corollary 1.5} Let $p$ be a prime.

{\rm (i)} For every $n=2,\ldots,p$ we have
$$\sum_{k=1}^n(-1)^{pk-1}\bi{pn-1}{pk-1}\eq(n-1)!B_{p-n}p^n\ (\mo\ p^{n+1}).\tag1.17$$

{\rm (ii)} Suppose that $p>3$. Then $p$ does not divide the class number $h_p$ of the
$p$-th cyclotomic field
$\Q(\zeta_p)$, if and only if
$$\ord_p\(\sum_{k=1}^n(-1)^k\bi{pn-1}{pk-1}\)=n\ \ \t{for all}\ n=3,5,\ldots,p-2.$$
Also,
$$\aligned&\sum_{k=1}^{(p-1)/2}(-1)^{k-1}\bi{p(p-1)/2-1}{pk-1}
\\\eq&[\![4\mid p+1]\!](-1)^{(h(-p)+1)/2}h(-p)p^{(p-1)/2}\ \l(\mo\ p^{(p+1)/2}\r),
\endaligned\tag1.18$$
where $h(-p)$ is the class number of the imaginary quadratic field
$\Q(\sqrt{-p})$.
\endproclaim
\Proof. (i) Let $n\in\{2,\ldots,p\}$. Then $\lfloor(pn-1-1)/(p-1)\rfloor=n$
and hence
$$F_p(pn-1,-1)=(-p)^{-n}C_p(pn-1,-1)=(-p)^{-n}\sum_{k=1}^n\bi{pn-1}{pk-1}(-1)^{pk-1}.$$
By Corollary 1.4,
$(-1)^nF_p(pn-1,-1)$
is congruent to
$$(p-1-\{-(n-1)\}_{p-1})!B_{\{-(n-1)\}_{p-1}}(-1+1)
=(n-1)!B_{p-n}$$
modulo $p$. Therefore (1.17) holds.

(ii) In view of part (i),
$$\align &\ord_p\(\sum_{k=1}^n(-1)^k\bi{pn-1}{pk-1}\)=n\ \ \t{for}\ n=3,5,\ldots,p-2
\\\iff&B_{p-n}\not\eq0\ (\mo\ p)\ \ \t{for}\ n=3,5,\ldots,p-2
\\\iff&p\ \t{is regular}
\\\iff& h_p\not\eq0\ (\mo\ p).
\endalign$$

 Taking $n=(p-1)/2$ in (1.17) we get
$$\align&\sum_{k=1}^{(p-1)/2}(-1)^{k-1}\bi{p(p-1)/2-1}{pk-1}
\\\eq&\f{((p-1)/2)!}{(p-1)/2}p^{(p-1)/2}B_{(p+1)/2}\ \l(\mo\ p^{(p+1)/2}\r).
\endalign$$
If $p\eq1\ (\mo\ 4)$, then $B_{(p+1)/2}=0$ since $(p+1)/2\in\{3,5,\ldots\}$.
If $p\eq 3\ (\mo\ 4)$, then $h(-p)\eq-2B_{(p+1)/2}\ (\mo\ p)$ (cf. [IR, p.\,238]),
and $((p-1)/2)!\eq(-1)^{(h(-p)+1)/2}\ (\mo\ p)$ by Mordell [M].
So (1.18) follows from the above.
This concludes the proof. \qed

\Remark. Let $p$ be an odd prime. If $p\gs5$, then (1.17) in the case $n=2$
reduces to Wolstenholme's congruence $\bi{2p-1}{p-1}\eq 1\ (\mo\ p^3)$ since $B_{p-2}=0$.
Taking $n=3$ in (1.17) we get
$$\bi{3p-1}{p-1}-\bi{3p-1}{2p-1}+\bi{3p-1}{3p-1}\eq 2B_{p-3}p^3\ (\mo\ p^4);$$
as $\bi{3p-1}{2p-1}=2\bi{3p-1}{p-1}$ this yields the congruence
$$\bi{3p-1}{p-1}\eq1-2p^3B_{p-3}\ (\mo\ p^4).$$
This was first obtained by
J.W.L. Glaisher (cf. [G1, p.\,21] and [G2, p.\,323]) who showed that
$$\bi{pn-1}{p-1}\eq1-\f{n(n-1)}3p^3B_{p-3}\ (\mo\ p^4)\ \ \t{for}\ n=1,2,3,\ldots.$$

\proclaim{Corollary 1.6} Let $p$ be an odd prime, and let
$n\in\{3,\ldots,p\}$ and $r\in\Z$. Then
$$F_p(pn-2,r)\eq-n!\l(\f{B_{p-n+1}(-r)}{n-1}+(r+1)\f{B_{p-n}(-r)}n\r)\ (\mo\ p).\tag1.19$$
\endproclaim
\Proof.  Clearly $\{-(pn-2)\}_{p-1}=p-n+1$. By Theorem 1.2, $F_p(pn-2,r)$
is congruent to
$$-(p-1-(p-n+1))!B_{p-n+1}^{(2)}(-r)
=-(n-2)!B_{p-n+1}^{(2)}(-r)$$
modulo $p$.

Let $m=p-n+1$. By [PS, (2.14)] or [SP, (1.12)],
$$\align&\f{(-1)^m}m\sum_{k=0}^m\bi mkB_kB_{m-k}(x)-\f{B_m(1-x)}{m}B_0
\\=&-\sum_{k=0}^1\bi 1kB_{1-k}(x)B_{m-1+k}(1-x)-B_1B_{m-1}(1-x)
\\=&-B_1(x)B_{m-1}(1-x)-B_0(x)B_m(1-x)-B_1B_{m-1}(1-x)
\\=&(-1)^m\l((B_1(x)+B_1)B_{m-1}(x)-B_m(x)\r)
\\=&(-1)^m\l((x-1)B_{m-1}(x)-B_m(x)\r).
\endalign$$
It follows that
$$\align B^{(2)}_{m}(-r)=&\sum_{k=0}^{m}\bi mk B_kB_{m-k}(-r)
\\=&(1-m)B_m(-r)+m(-r-1)B_{m-1}(-r)
\\\eq&(1+n-1)B_{p-n+1}(-r)-(r+1)(-n+1)B_{p-n}(-r)
\\\eq&n(n-1)\l(\f{B_{p-n+1}(-r)}{n-1}+(r+1)\f{B_{p-n}(-r)}n\r)\ (\mo\ p).
\endalign$$

Combining the above we immediately obtain (1.19). \qed

\medskip

By Theorem 1.1 or 1.2, for any prime $p$ the Fleck quotient
$F_p(n,r)$ (with $n\in\N$ and $r\in\Z$) modulo $p$
only depends on $p$ and $r$ and the remainder of $n$ modulo $p(p-1)$.
This observation can be further extended as follows.

\proclaim{Theorem 1.3} Let $p$ be a prime, and let $a,l,n\in\N$ and $r\in\Z$.
Then
$$\aligned&\sum_{k=0}^n\bi nk(-1)^kF_p\l(kp^a(p-1)+l,r\r)
\\&\eq0
\ \l(\mo\ p^{an+\lceil(n-l^*)/(p-1)\rceil}\r),
\endaligned\tag1.20$$
where $l^*=\{-l\}_{p-1}$ and $\lceil\cdot\rceil$ is the ceiling function.
\endproclaim

The following consequence is somewhat similar to Kummer's congruence
for Bernoulli numbers (cf. [IR, pp.\,238--241]).

\proclaim{Corollary 1.7} Let $p$ be a prime, and let $a,l\in\N$ and $r\in\Z$.
Then
$$\align F_p(p^a(p-1)+l,r)\eq& F_p(l,r)\ (\mo\ p^a),
\\F_p(2p^a(p-1)+l,r)\eq&2F_p(p^a(p-1)+l,r)-F_p(l,r)\ (\mo\ p^{2a}),
\\F_p(3p^a(p-1)+l,r)\eq&3F_p(2p^a(p-1)+l,r)-3F_p(p^a(p-1)+l,r)
\\&+F_p(l,r)\ \ (\mo\ p^{3a}).
\endalign$$
\endproclaim
\Proof. Simply apply (1.20) with $n=1,2,3$. \qed

\medskip
Let $p$  be a prime, and let $a\in\Z^+$ and $r\in\Z$. In 1977
C. S. Weisman [We] extended Fleck's result by showing that
if $n\gs p^{a-1}$ then
$$C_{p^a}(n,r)\eq0\ \l(\mo\ p^{\lfloor(n-p^{a-1})/\varphi(p^a)\rfloor}\r),$$
where $\varphi$ is Euler's totient function.
In view of this, we define the {\it generalized Fleck quotient}
$$F_{p^a}(n,r)=(-p)^{-\lfloor(n-p^{a-1})/\varphi(p^a)\rfloor}C_{p^a}(n,r)
+[\![n<p^{a-1}]\!]\in\Z.$$
Note that $F_{p^a}(n,r)\eq1\ (\mo\ p)$ for $n=0,\ldots,p^{a-1}-1$.

\proclaim{Theorem 1.4} Let $p$ be a prime, and let $a,n\in\Z^+$ with $n\gs p^{a-1}$.

{\rm (i)} For any $r\in\Z$ we have
$$F_{p^a}(n,r)\eq\sum_{k=0}^{d}\bi{r+k-1}kF_{p^a}(n+k,0)\ (\mo\ p),\tag1.21$$
where $d=\{p^{a-1}-1-n\}_{\varphi(p^a)}$ is the least nonnegative integer
with $n+d\eq p^{a-1}-1\ (\mo\ \varphi(p^a))$.

{\rm (ii)} We have
$$\ord_p\l(C_{p^a}(n,r)\r)=\l\lfloor\f{n-p^{a-1}}{\varphi(p^a)}\r\rfloor
\ (\t{i.e.,}\ p\nmid F_{p^a}(n,r))\quad\t{for some}\ r\in\Z.\tag1.22$$
If $n\gs 2p^{a-1}$, then
$$F_{p^a}\l(n+p^a(p-1),r\r)\eq F_{p^a}(n,r)\ (\mo\ p)\quad\t{for all}\ r\in\Z.\tag1.23$$
\endproclaim

In view of the first congruence in Corollary 1.7 and
the last congruence in Theorem 1.4, we propose the following conjecture.
\proclaim{Conjecture 1.1} Let $p$ be a prime, and let $a,b,n\in\Z^+$ and $r\in\Z$.
If $n\gs 2p^{a+b-2}$, then
$$F_{p^a}\l(n+\varphi(p^{a+b}),r\r)\eq F_{p^a}(n,r)\ (\mo\ p^b).$$
\endproclaim

Theorems 1.1, 1.2 and 1.3 will be proved in Sections 2, 3 and 4 respectively.
In Section 5 we will first give a new proof of Weisman's congruence via roots of unity,
and then establish Theorem 1.4.

\heading{2. Proof of Theorem 1.1}\endheading

\proclaim{Lemma 2.1} Let $p$ be a prime, and let $n\in\N$ and $n^*=\{-n\}_{p-1}$.
Define $G(n)=\sum_{a=1}^{p-1}a^n\zeta_p^a$ and $\pi=1-\zeta_p$,
where $\zeta_p$ is a primitive $p$-th root of unity in the complex field $\C$.
Then
$$G(n)\eq(-1)^{n^*-1}\sum_{m=n^*}^{p-2}s(m,n^*)\f{\pi^m}{m!}\ (\mo\ p),\tag2.1$$
where $s(m,0),\ldots,s(m,m)$ are Stirling numbers of the first kind
defined by $(x)_m=\sum_{k=0}^m(-1)^{m-k}s(m,k)x^k$.
\endproclaim
\Proof. Clearly,
$$\align G(n)=&\sum_{a=1}^{p-1}a^n(1-\pi)^a=\sum_{a=1}^{p-1}a^n\sum_{m=0}^a\bi am(-\pi)^m
\\=&\sum_{m=0}^{p-1}\f{(-\pi)^m}{m!}\sum_{a=1}^{p-1}a^n(a)_m
\\=&\sum_{m=0}^{p-1}\f{(-\pi)^m}{m!}\sum_{a=1}^{p-1}a^n\sum_{k=0}^m(-1)^{m-k}s(m,k)a^k
\\=&\sum_{m=0}^{p-1}\f{(-\pi)^m}{m!}\sum_{k=0}^m(-1)^{m-k}s(m,k)\sum_{a=1}^{p-1}a^{n+k}.
\endalign$$

Since
$$1+x+\cdots+x^{p-1}=\f{x^p-1}{x-1}=\prod_{a=1}^{p-1}(x-\zeta_p^a),$$
we have
$$\f{p}{\pi^{p-1}}=\prod_{a=1}^{p-1}\f{1-\zeta_p^a}{\pi}
=\prod_{a=1}^{p-1}\f{1-(1-\pi)^a}{\pi}
\eq\prod_{a=1}^{p-1}a\eq-1\ (\mo\ \pi)$$
with the help of Wilson's theorem. Note also that
$$\sum_{a=1}^{p-1}a^{n+k}
\eq-[\![p-1\mid n+k]\!]\ (\mo\ p)$$
by elementary number theory (see, e.g., [IR, pp.\,235--236]).
Therefore
$$\align G(n)\eq&\sum_{m=0}^{p-2}\f{\pi^m}{m!}\sum_{k=0}^m(-1)^ks(m,k)(-[\![k=n^*]\!])
\\\eq&(-1)^{n^*-1}\sum_{m=n^*}^{p-2}s(m,n^*)\f{\pi^m}{m!}\ (\mo\ p).
\endalign$$
This concludes the proof. \qed

\Remark. Let $p$ be an odd prime. For each $a\in\Z$ let $\bar a=a+p\Z\in \Bbb F_p=\Z/p\Z$.
Let $\omega$ be the Teichm\"uller character of the multiplicative group $\Bbb F_p^*=\Bbb F_p\sm\{\bar0\}$.
For $\bar a\in \Bbb F_p^*$, $\omega(\bar a)$ is just the $(p-1)$-th root of unity
in the unique unramified extension of the $p$-adic field $\Q_p$
with $\omega(\bar a)\eq a\ (\mo\ p)$. (See, e.g., [Wa, p.\,51].)
If $\zeta_p$ is a primitive $p$-th root of unity in the algebraic closure
of $\Q_p$, then for $n\in\N$ and $\pi=1-\zeta_p$ we have
$$\sum_{a=1}^{p-1}a^n\zeta_p^a\eq\sum_{a=1}^{p-1}\omega^{n}(\bar a)\zeta_p^a
\eq-\f{(-\pi)^{n^*}}{n^*!}\ \ (\mo\ \pi^{n^*+1})$$
with $n^*=\{-n\}_{p-1}$, by Stickelberger's congruence for Gauss' sums (cf. [BEW, pp.\,344--345]).

\proclaim{Lemma 2.2} Let $p$ be a prime, and let $\zeta_p$ be a primitive $p$-th root of unity
in $\C$. Let $n=p^am+n_0>0$ with $a\in\Z^+$ and $m,n_0\in\N$.
Then, for any $r\in\Z$ we have
$$\align&\pi^{-p^am}C_p(n,r)-[\![p-1\mid m]\!]C_p(n_0,r)
\\\eq&\f{G(p^am)}p\sum_{k=0}^{n_0}\bi{n_0}k(-1)^k
(k-r)^{p^am^*}
\ \l(\mo\ p^{a-1}\pi^{\min\{n_0+1,\,p-1\}}\r),
\endalign$$
where $\pi=1-\zeta_p$ and $m^*=\{-m\}_{p-1}$.
\endproclaim
\Proof. Let $j\in\{1,\ldots,p-1\}$. Then
$$\l(\f{1-\zeta_p^j}{\pi}\r)^{m}=\l(\f{1-(1-\pi)^j}{\pi}\r)^{m}
=\(\sum_{i=1}^j\bi ji(-\pi)^{i-1}\)^{m}=j^{m}+\beta_j\pi,$$
where $\beta_j$ is a suitable element in the ring $\overline\Z$ of algebraic integers.
For $i=0,1,\ldots$, if
$$\l(\f{1-\zeta_p^j}{\pi}\r)^{p^im}=j^{p^im}+p^i\pi\beta_j^{(i)}$$
for some $\beta_j^{(i)}\in\overline\Z$, then
$$\l(\f{1-\zeta_p^j}{\pi}\r)^{p^{i+1}m}=\l(j^{p^im}+p^i\pi\beta_j^{(i)}\r)^p
=j^{p^{i+1}m}+p^{i+1}\pi\beta_j^{(i+1)}$$
for some $\beta_j^{(i+1)}\in\overline\Z$.
So
$$\l(\f{1-\zeta_p^j}{\pi}\r)^{p^am}\eq j^{p^am}\ (\mo\ p^a\pi).$$

Observe that
$$pC_p(n,r)=\sum_{j=0}^{p-1}\zeta_p^{-jr}(1-\zeta_p^j)^n
=\pi^{p^am}\sum_{j=1}^{p-1}\zeta_p^{-jr}
\l(\f{1-\zeta_p^j}{\pi}\r)^{p^am}(1-\zeta_p^j)^{n_0}.$$
As $\pi^{n_0}$ divides $(1-\zeta_p^j)^{n_0}$ in the ring $\overline\Z$,
by the above
$\pi^{-p^am}pC_p(n,r)$ is congruent to
$$\sum_{j=1}^{p-1}\zeta_p^{-jr}j^{p^am}
\sum_{k=0}^{n_0}\bi {n_0}k(-1)^k\zeta_p^{jk}
=\sum_{k=0}^{n_0}\bi{n_0}k(-1)^kS_{k-r}$$
modulo $p^a\pi^{n_0+1}$, where
$$S_{k-r}=\sum_{j=1}^{p-1}j^{p^am}\zeta_p^{j(k-r)}.$$

If $k\not\eq r\ (\mo\ p)$, then
$$\align S_{k-r}=&(k-r)^{-p^am}\sum_{j=1}^{p-1}(j(k-r))^{p^am}\zeta_p^{j(k-r)}
\\\eq&(k-r)^{p^am^*} \sum_{t=1}^{p-1}t^{p^am}\zeta_p^t
=(k-r)^{p^am^*} G(p^am)\ (\mo\ p^{a+1}).
\endalign$$
(Note that if $j(k-r)\eq t\ (\mo\ p)$ then  $(j(k-r))^{p^a}\eq t^{p^a}\ (\mo\ p^{a+1})$.)

Choose a primitive root $g$ modulo $p$. Since
$$(g^{p^am}-1)\sum_{j=1}^{p-1}j^{p^am}=\sum_{j=1}^{p-1}(gj)^{p^am}
-\sum_{t=1}^{p-1}t^{p^am}\eq0\ (\mo\ p^{a+1}),$$
if $p-1\nmid m$ then $g^{p^am}-1\not\eq0\ (\mo\ p)$ and so
$\sum_{j=1}^{p-1}j^{p^am}\eq0\ (\mo\ p^{a+1})$.
Thus, when $k\eq r\ (\mo\ p)$ we have
$$S_{k-r}=\sum_{j=1}^{p-1}j^{p^am}\eq (p-1)[\![p-1\mid m]\!]
\ (\mo\ p^{a+1}).$$

Recall that $p/\pi^{p-1}\eq-1\ (\mo\ \pi)$. In view of the above,
$$\align&\pi^{-p^am}pC_p(n,r)-\sum_{k=0}^{n_0}\bi {n_0}k(-1)^k(k-r)^{p^am^*}G(p^am)
\\\eq&\sum^{n_0}\Sb k=0\\p\mid k-r\endSb\bi{n_0}k(-1)^k
\l([\![p-1\mid m]\!](p-1)-(k-r)^{p^am^*}G(p^am)\r)
\\\eq&C_p(n_0,r)[\![p-1\mid m]\!]p\ \ \l(\mo\ p^a\pi^{\min\{n_0+1,\,p-1\}}\r),
\endalign$$
where we have noted that if $p-1\mid m$ (i.e., $m^*=0$) then
$$p-1-G(p^am)\eq p-\sum_{t=0}^{p-1}\zeta_p^t
=p-\f{1-\zeta_p^p}{1-\zeta_p}=p\ (\mo\ p^{a+1}).$$
Therefore the desired congruence follows. \qed

\medskip
\noindent{\it Proof of Theorem 1.1}. In the case $n=0$, (1.2) holds
since $n_1=n_0=0$ and $F_p(n,r)=-pC_p(0,r)+1$.
Below we assume $n>0$.

Let $\zeta_p$ be a primitive $p$-th root of unity in $\C$,
and set $\pi=1-\zeta_p$. By Lemma 2.2 in the case $a=1$,
$$\align&\pi^{-p\lfloor n/p\rfloor}C_p(n,r)-[\![n_1=0]\!]C_p(n_0,r)
\\\eq&\f{G(p\lfloor n/p\rfloor)}p\sum_{k=0}^{n_0}\bi{n_0}k(-1)^k(k-r)^{pn_1}
\ \l(\mo\ \pi^{\min\{n_0+1,\,p-1\}}\r).
\endalign$$
In view of Lemma 2.1,
$$G\l(p\l\lfloor\f np\r\rfloor\r)
\eq G\l(\l\lfloor \f np\r\rfloor\r)
\eq(-1)^{n_1-1}\sum_{m=n_1}^{p-2}s(m,n_1)\f{\pi^m}{m!}\ (\mo\ p).$$
If $n_0>n_1$, then
$$\sum_{k=0}^{n_0}\bi{n_0}k(-1)^k(k-r)^{pn_1}
\eq\sum_{k=0}^{n_0}\bi{n_0}k(-1)^k(k-r)^{n_1}=0\ (\mo\ p),$$
where we have applied Fermat's little theorem and Euler's identity (mentioned in Section 1).
Therefore
$$\align&\pi^{-p\lfloor n/p\rfloor}C_p(n,r)-[\![n_1=0]\!]C_p(n_0,r)
\\\eq&\f{(-1)^{n_1-1}}p\sum_{m=n_1}^{p-2}s(m,n_1)\f{\pi^m}{m!}
\sum_{k=0}^{n_0}\bi{n_0}k(-1)^k(k-r)^{pn_1}
\\&\qquad\qquad\qquad\qquad\ \l(\mo\ \pi^{[\![n_0>n_1]\!]\min\{n_0+1,\,p-1\}}\r).
\endalign$$
Recall that $-p/\pi^{p-1}\eq1\ (\mo\ \pi)$.
Since $s(n_1,n_1)=1$ and
$$\f{p^{[\![n_0\ls n_1]\!]}}{\pi^{n_1}}\pi^{[\![n_0>n_1]\!]\min\{n_0+1,\,p-1\}}
\eq0\ (\mo\ \pi),$$
by the above we have
$$\align&\f{p^{[\![n_0\ls n_1]\!]}C_p(n,r)}{\pi^{p\lfloor n/p\rfloor+n_1}}
-p^{[\![n_0=0]\!]}[\![n_1=0]\!]C_p(n_0,r)
\\\eq&\f{(-1)^{n_1-1}/n_1!}{p^{[\![n_0>n_1]\!]}}\sum_{k=0}^{n_0}\bi{n_0}k(-1)^k(k-r)^{pn_1}
\ (\mo\ \pi).
\endalign$$

 Note that
$$\l\lfloor\f{n-1}{p-1}\r\rfloor=\l\lfloor\f{p\lfloor n/p\rfloor+n_0-1}{p-1}\r\rfloor
=\f{p\lfloor n/p\rfloor+n_1}{p-1}-[\![n_0\ls n_1]\!]$$
and hence
$$\align\f{(-p)^{[\![n_0\ls n_1]\!]}C_p(n,r)}{\pi^{p\lfloor n/p\rfloor+n_1}}
=&\f{C_p(n,r)}{(-p)^{\lfloor(n-1)/(p-1)\rfloor}}
\l(\f{-p}{\pi^{p-1}}\r)^{(p\lfloor n/p\rfloor+n_1)/(p-1)}
\\\eq& F_p(n,r)\ (\mo\ \pi).
\endalign$$
In view of the above,
$$\align&(-1)^{[n_0\ls n_1]}F_p(n,r)-[\![n_0>n_1=0]\!]C_p(n_0,r)
\\\eq&\f{(-1)^{n_1-1}/n_1!}{p^{[\![n_0>n_1]\!]}}\sum_{k=0}^{n_0}\bi{n_0}k(-1)^k(k-r)^{pn_1}
\ (\mo\ \pi).
\endalign$$
As the rational $p$-adic integer
$$\align D=&F_p(n,r)-[\![n_0>n_1=0]\!]C_p(n_0,r)
\\&-\f{(-1)^{n_1}}{(-p)^{[\![n_0>n_1]\!]}\cdot n_1!}\sum_{k=0}^{n_0}\bi{n_0}k(-1)^k(k-r)^{pn_1}
\endalign$$
is divisible by $\pi$, we have $D^{p-1}\eq0\ (\mo\ p)$ and hence $D\eq0\ (\mo\ p)$.
Thus
$$\aligned&F_p(n,r)-[\![n_0>n_1=0]\!]C_p(n_0,r)
\\\eq&\f{(-1)^{n_1}}{(-p)^{[\![n_0>n_1]\!]}\cdot n_1!}
\sum_{k=0}^{n_0}\bi{n_0}k(-1)^k(k-r)^{pn_1}
\ (\mo\ p).
\endaligned\tag2.2$$

In the case $n_0\ls n_1$, (2.2) reduces to (1.2).
When $n_0>n_1=0$, (2.2) yields (1.3) since
$C_p(n_0,r)=(-1)^{\{r\}_p}\bi{n_0}{\{r\}_p}$
and $\sum_{k=0}^{n_0}\bi{n_0}k(-1)^k=(1-1)^{n_0}=0$.

 Now assume that $n_0>n_1>0$. As $\sum_{k=0}^{n_0}\bi{n_0}k(k-r)^{n_1}=0$
 by Euler's identity, (2.2) implies that
 $$F_p(n,r)\eq\f{(-1)^{n_1-1}}{n_1!}
 \sum_{k=0}^{n_0}\bi{n_0}k(-1)^{k}\f{(k-r)^{pn_1}-(k-r)^{n_1}}p\ (\mo\ p).$$
If $n_1=1$, then
$$\f{(k-r)^{pn_1}-(k-r)^{n_1}}p=(k-r)^{n_1}n_1q_p(k-r);$$
if $n_1\gs2$ and $k\eq r\ (\mo\ p)$, then
$$\f{(k-r)^{pn_1}-(k-r)^{n_1}}p\eq0\eq(k-r)^{n_1}n_1q_p(k-r)\ (\mo\ p);$$
if $a=k-r\not\eq0\ (\mo\ p)$, then
$$\f{(k-r)^{pn_1}-(k-r)^{n_1}}p
=a^{n_1}\f{(1+p\cdot q_p(a))^{n_1}-1}p\eq a^{n_1}n_1q_p(a)\ (\mo\ p).$$
Therefore (1.4) follows.

The proof is now complete. \qed

\heading{3. Proof of Theorem 1.2}\endheading

The following lemma is a refinement of
an induction technique used by Sun [S06].

\proclaim{Lemma 3.1} Let $p$ be a prime, and let $n\in\N$ with $n\gs p$. Then
$$F_p(n,r)\eq-\sum_{j=1}^{p-1}\f1j\sum_{i=0}^{j-1}F_p(n-p+1,r-i)\ (\mo\ p).\tag3.1$$
\endproclaim
\Proof. Set $n'=n-(p-1)>0$. By the Chu-Vandermonde convolution identity (cf. [GKP, (5.27)]),
$$\align F_p(n,r)=&(-p)^{-\lfloor(n-1)/(p-1)\rfloor}\sum\Sb 0\ls k\ls n\\k\eq r\,(\mo\ p)\endSb
\sum_{j=0}^{k}\bi{p-1}j\bi{n'}{k-j}(-1)^k
\\=&-\f1p\sum_{j=0}^{p-1}\bi{p-1}j
(-p)^{-\lfloor(n'-1)/(p-1)\rfloor}\sum\Sb j\ls k\ls n\\p\mid k-r\endSb\bi {n'}{k-j}(-1)^{k}
\\=&-\f1p\sum_{j=0}^{p-1}\bi{p-1}j(-1)^jF_p(n',r-j).
\endalign$$

For any $j=0,\ldots,p-1$,
clearly
$$\align&\bi{p-1}j(-1)^j=\prod_{0<i\ls j}\l(1-\f pi\r)
\\\eq&1-\sum_{0<i\ls j}\f pi\eq(-1)^{p-1}+p\sum_{j<k<p}\f 1k\ (\mo\ p^2).
\endalign$$
(Note that $2\sum_{k=1}^{p-1}1/k=\sum_{k=1}^{p-1}(1/k+1/(p-k))\eq0\ (\mo\ p)$.)
Also,
$$\sum_{j=0}^{p-1}F_p(n',r-j)=(-p)^{-\lfloor(n'-1)/(p-1)\rfloor}
\sum_{k=0}^{n'}\bi{n'}k(-1)^k=0.$$
Therefore
$$F_p(n,r)\eq-\sum_{j=0}^{p-1}\sum_{j<k<p}\f{F_p(n',r-j)}k
=-\sum_{k=1}^{p-1}\f1k\sum_{j=0}^{k-1}F_p(n',r-j)\ (\mo\ p).$$
This proves (3.1). \qed

\medskip
\noindent{\it Proof of Theorem 1.2}.
(i) Suppose $m\gs0$. Then
$$\align&\sum_{k=0}^{n^*}\bar S(m+n^*-k,m)\f{(-r)^k}{k!}
\\=&[x^{m+n^*}]\sum_{l=m}^{\infty}\bar S(l,m)x^l\sum_{k=0}^{\infty}\f{(-rx)^k}{k!}
\\=&[x^{m+n^*}](e^x-1)^me^{-rx}=[x^{n^*}]\l(\f{e^x-1}x\r)^me^{-rx}
\\=&[x^{m+n^*}]\sum_{k=0}^m\bi mk(-1)^{m-k}e^{(k-r)x}
=\sum_{k=0}^m\bi mk(-1)^{m-k}\f{(k-r)^{m+n^*}}{(m+n^*)!}.
\endalign$$
By the identity (2.4) of Sun [S03], for any $l=0,1,\ldots$ we have
$$\align\sum_{k=0}^{m}\bi{m}k(-1)^{m-k}(k+l)^{m+n^*}
=&\sum_{j=0}^{l}\bi{l}j(m+j)!S(m+n^*,m+j)
\\=&\sum_{j=0}^{n^*}\bi{l}j(m+j)!S(m+n^*,m+j).
\endalign$$
Thus
$$\sum_{k=0}^{m}\bi{m}k(-1)^{m-k}(k+x)^{m+n^*}=\sum_{j=0}^{n^*}\bi{x}j(m+j)!S(m+n^*,m+j)$$
and hence
$$\sum_{k=0}^{m}\bi{m}k(-1)^{m-k}\f{(k-r)^{m+n^*}}{(m+n^*)!}
=\sum_{j=0}^{n^*}\bi{-r}j\bar S(m+n^*,m+j).$$

If $m\ls 0$, then
$$\f{B_{n^*}^{(-m)}(-r)}{n^*!}=[x^{n^*}]\l(\f x{e^x-1}\r)^{-m}e^{-rx}
=[x^{n^*}]\l(\f {e^x-1}x\r)^{m}e^{-rx}.$$
Note also that
$$\f{1}{n^*!}=\f{\prod_{j=1}^{p-1-n^*}(p-j)}{(p-1)!}
\eq(-1)^{n^*+1}(p-1-n^*)!\ (\mo\ p)$$
by Wilson's theorem.

In view of the above, whether $m\gs0$ or $m\ls 0$, we only need to show that
$$(-1)^nF_p(n,r)\eq[x^{n^*}]\l(\f{e^x-1}x\r)^{m}e^{-rx}\ (\mo\ p).$$

(ii) All those formal power series $f(x)=\sum_{k=0}^{\infty}a_kx^k$ with $a_k\in\Q$
and $a_0,\ldots,a_{n^*}\in\Z_p$ form a ring $R_{n^*}$
under the usual addition and multiplication.
In particular, this ring contains
$$e^{-rx}=\sum_{k=0}^{\infty}(-r)^k\f{x^k}{k!},\ \ \f{e^x-1}x=\sum_{k=0}^{\infty}\f{x^k}{(k+1)!}
\ \ \t{and}\ \ \f x{e^x-1}=\sum_{k=0}^{\infty}B_k\f{x^k}{k!}.$$
(Recall that $n^*<p-1$ and $B_0,\ldots,B_{n^*}\in\Z_p$.)
If $f(x)=\sum_{k=0}^{\infty}a_kx^k$ and $g(x)=\sum_{k=0}^{\infty}b_kx^k$
belong to $R_{n^*}$, then
$$\align&[x^{n^*}]f(x)g(x)^p
=[x^{n^*}]\sum_{j=0}^{n^*}a_jx^j\l(\sum_{k=0}^{n^*}b_kx^k\r)^p
\\\eq&[x^{n^*}]\sum_{j=0}^{n^*}a_jx^j
\sum_{k=0}^{n^*}b_k^px^{pk}
=a_{n^*}b_0^p\eq[x^{n^*}]f(x)[x^0]g(x)\ (\mo\ p).
\endalign$$
Consequently, for any $a\in\Z$ we have
$$[x^{n^*}]\l(\f{e^x-1}x\r)^me^{ax}\eq[x^{n^*}]\l(\f{e^x-1}x\r)^ne^{ax}\ (\mo\ p)$$
since $m\eq n\ (\mo\ p)$.
By this and part (i), it suffices to use induction on $n$ to show that
$$(-1)^nF_p(n,r)\eq[x^{n^*}]\l(\f{e^x-1}x\r)^{n}e^{-rx}\ (\mo\ p).\tag3.2$$

(iii) Obviously
$$(-1)^0F_p(0,r)=-pC_p(0,r)+1\eq1=[x^0]\l(\f{e^x-1}x\r)^0e^{-rx}\ (\mo\ p).$$
So (3.2) holds for $n=0$.

Suppose that $0<n\ls p-1$. Then $n^*=p-1-n$ and
$$\align&[x^{n^*}]\l(\f{e^x-1}x\r)^ne^{-rx}
=[x^{p-1}](e^x-1)^ne^{-rx}
\\=&\sum_{k=0}^n\bi nk(-1)^{n-k}[x^{p-1}]e^{(k-r)x}
=\sum_{k=0}^n\bi nk(-1)^{n-k}\f{(k-r)^{p-1}}{(p-1)!}
\\\eq&(-1)^{n-1}\sum_{k\not\eq r\,(\mo\ p)}\bi nk(-1)^k\ (\mo\ p).
\endalign$$
(To get the last congruence we have applied Wilson's theorem and Fermat's little theorem.)
Since
$$-\sum_{k\not\eq r\,(\mo\ p)}\bi nk(-1)^k
=\sum_{k\eq r\,(\mo\ p)}\bi nk(-1)^k=F_p(n,r),$$
the desired (3.2) follows.

Now fix $n\gs p$ and assume that (3.2) holds for smaller values of $n$.
Clearly $n'=n-(p-1)>0$ and $\{-n'\}_{p-1}=n^*$. In light of Lemma 3.1,
$$F_p(n,r)\eq-\sum_{j=1}^{p-1}\f1j\sum_{k=0}^{j-1}F_p(n',r-k)\ (\mo\ p).$$
By the induction hypothesis and part (ii),
$$\align(-1)^{n'}F_p(n',r-k)\eq&[x^{n^*}]\l(\f {e^x-1}x\r)^{n'}e^{-(r-k)x}
\\\eq&[x^{n^*}]\l(\f {e^x-1}x\r)^{n+1}e^{(k-r)x}\ (\mo\ p).
\endalign$$
Thus $(-1)^{n-1}F_p(n,r)$
is congruent to
$$\align&\sum_{j=1}^{p-1}\f1j\sum_{k=0}^{j-1}\([x^{n^*}]
\l(\f {e^x-1}x\r)^{n+1}e^{(k-r)x}\)
\\=&[x^{n^*}]\l(\f {e^x-1}x\r)^{n+1}e^{-rx}\sum_{j=1}^{p-1}\(\f1j\cdot\f{e^{jx}-1}{e^x-1}\)
\\=&[x^{n^*}]\l(\f {e^x-1}x\r)^{n}e^{-rx}\sum_{j=1}^{p-1}\f{e^{jx}-1}{jx}
\endalign$$
modulo $p$. This yields
$$\align (-1)^nF_p(n,r)\eq&-[x^{n^*}]\l(\f {e^x-1}x\r)^{n}e^{-rx}
\sum_{j=1}^{p-1}\sum_{k=1}^{p-1}\f{(jx)^{k-1}}{k!}
\\\eq&[x^{n^*}]\l(\f {e^x-1}x\r)^{n}e^{-rx}\ (\mo\ p),
\endalign$$
since $n^*<p-1$ and
$\sum_{j=1}^{p-1}j^{k-1}\eq-[\![p-1\mid k-1]\!]\ (\mo\ p)$.
\medskip

In view of the above, we have completed the proof. \qed

\heading{4. Proof of Theorem 1.3}\endheading

\noindent{\it Proof of Theorem 1.3}.
Let $\zeta_p$ be a primitive $p$-th root of unity in $\C$, and set $\pi=1-\zeta_p$.
For any $k=0,\ldots,n$, we have
$$\align pC_p(kp^a(p-1)+l,r)=&\sum_{j=0}^{p-1}\zeta_p^{-jr}(1-\zeta_p^j)^{kp^a(p-1)+l}
\\=&\sum_{j=1}^{p-1}\zeta_p^{-jr}(1-\zeta_p^j)^{kp^a(p-1)+l}+[\![k=l=0]\!]
\endalign$$
and thus
$$\align &F_p(kp^a(p-1)+l,r)
\\=&(-p)^{-\lfloor(kp^a(p-1)+l-1)/(p-1)\rfloor}C_p(kp^a(p-1)+l,r)+[\![k=l=0]\!]
\\=&-(-p)^{-kp^a-\lfloor(l-1)/(p-1)\rfloor-1}
\sum_{j=1}^{p-1}\zeta_p^{-jr}(1-\zeta_p^j)^{kp^a(p-1)+l}.
\endalign$$
Therefore, for $S_n=\sum_{k=0}^n\bi nk(-1)^kF_p(kp^a(p-1)+l,r)$ we have
$$S_n=-\sum_{j=1}^{p-1}\zeta_p^{-jr}(1-\zeta_p^j)^l
(-p)^{-\lfloor(l-1)/(p-1)\rfloor-1}c_{n,j},\tag4.1$$
where
$$\align c_{n,j}=&\sum_{k=0}^n\bi nk(-1)^k(-p)^{-kp^a}(1-\zeta_p^j)^{kp^a(p-1)}
\\=&\l(1-(-p)^{-p^a}(1-\zeta_p^j)^{p^a(p-1)}\r)^n.
\endalign$$

Let $j\in\{1,\ldots,p-1\}$. Clearly
$$\(\f{1-\zeta_p^j}{\pi}\)^{p-1}=\(\f{1-(1-\pi)^j}{\pi}\)^{p-1}\eq j^{p-1}\eq1\ (\mo\ \pi)$$
and hence
$$b_j:=\f{(1-\zeta_p^j)^{p-1}}{-p}
=\l(\f{1-\zeta_p^j}{\pi}\r)^{p-1}\f{\pi^{p-1}}{-p}\eq1\ (\mo\ \pi).$$
(Recall the congruence $p/\pi^{p-1}\eq-1\ (\mo\ \pi)$.) It follows that
$b_j^{p^a}\eq1\ (\mo\ p^a\pi)$ and
$$c_{n,j}=\l(1-b_j^{p^a}\r)^n\eq0\ (\mo\ p^{an}\pi^n).\tag4.2$$

Since $(1-\zeta_p^j)^l\eq0\ (\mo\ \pi^l)$ and $\ord_p(\pi)=1/(p-1)$,
in view of (4.1) and (4.2) we have
$$\ord_p(S_n)\gs\f{l+n}{p-1}+an-\l\lfloor\f{l-1}{p-1}\r\rfloor-1
=an+\f{l+n}{p-1}-\f{l+l^*}{p-1}=an+\f{n-l^*}{p-1}$$
and hence $\ord_p(S_n)\gs an+\lceil(n-l^*)/(p-1)\rceil$.
This proves (1.20). \qed

\heading{5. On generalized Fleck quotients}\endheading

\proclaim{Lemma 5.1}
Let $d,q\in\Z^+$, $n\in\N$ and $r\in\Z$.
Let $\zeta_{dq}$ be a primitive $dq$-th root of unity in $\C$. Then
$$C_{dq}(n,r)=\f1d\sum_{k=0}^n\bi nkC_q(k,r)
\sum_{j=0}^{d-1}\zeta_{dq}^{j(k-r)}\l(1-\zeta_{dq}^j\r)^{n-k}.\tag5.1$$
\endproclaim
\Proof. Note that $\zeta=\zeta_{dq}^d$ is a primitive $q$-th root of unity.
Thus
$$\align&q\sum_{k=0}^n\bi nkC_q(k,r)
\sum_{j=0}^{d-1}\zeta_{dq}^{j(k-r)}\l(1-\zeta_{dq}^j\r)^{n-k}
\\=&\sum_{k=0}^n\bi nk\sum_{s=0}^{q-1}\zeta^{-sr}(1-\zeta^s)^k
\sum_{j=0}^{d-1}\zeta_{dq}^{j(k-r)}\l(1-\zeta_{dq}^j\r)^{n-k}
\\=&\sum_{s=0}^{q-1}\sum_{j=0}^{d-1}\zeta_{dq}^{-(ds+j)r}\sum_{k=0}^n\bi nk
\l(\zeta_{dq}^j(1-\zeta_{dq}^{ds})\r)^k\l(1-\zeta_{dq}^j\r)^{n-k}
\\=&\sum_{s=0}^{q-1}\sum_{j=0}^{d-1}\zeta_{dq}^{-(ds+j)r}\l(1-\zeta_{dq}^{ds+j}\r)^n
\\=&\sum_{t=0}^{dq-1}\zeta_{dq}^{-tr}\l(1-\zeta_{dq}^t\r)^n=dqC_{dq}(n,r).
\endalign$$
So we have (5.1). \qed

With the help of Lemma 5.1 we can prove the following result
via roots of unity.

\proclaim{Theorem 5.1 {\rm (Weisman, 1977)}} Let $p$ be a prime, and let
$a\in\Z^+$, $n\in\N$ and $r\in\Z$. Then $F_{p^a}(n,r)\in\Z$.
\endproclaim
\Proof. We use induction on $a$.

The case $a=1$ reduces to Fleck's result.
A proof of Fleck's result via roots of unity was given by A. Granville [Gr].

Now let $a\gs 2$ and assume that $F_{p^{a-1}}(n',r')\in\Z$ for all $n'\in\N$ and $r'\in\Z$.
If $n<p^a$, then $\lfloor(n-p^{a-1})/\varphi(p^a)\rfloor\ls0$ and hence $F_{p^a}(n,r)\in\Z$.
Below we suppose $n\gs p^a$ and let $\zeta_{p^a}$ be a primitive $p^a$-th root of unity
in $\C$.

 By Lemma 5.1,
$$C_{p^a}(n,r)=\f1p\sum_{k=0}^n\bi nkC_{p^{a-1}}(k,r)
\sum_{j=0}^{p-1}\zeta_{p^a}^{j(k-r)}\l(1-\zeta_{p^a}^j\r)^{n-k}.\tag5.2$$
Observe that
$$\prod^{p^a-1}\Sb j=1\\p\nmid j\endSb\l(1-\zeta_{p^a}^j\r)
=\prod\Sb \gamma^{p^a}=1\\\gamma^{p^{a-1}}\not=1\endSb(1-\gamma)
=\lim_{x\to 1}\f{x^{p^a}-1}{x^{p^{a-1}}-1}=\f{p^a}{p^{a-1}}=p.$$
If $p\nmid j$, then $(1-\zeta_{p^a}^j)/(1-\zeta_{p^a})$
is a unit in the ring $\Z[\zeta_{p^a}]$
and thus
$$\ord_p(1-\zeta_{p^a}^j)=\ord_p(1-\zeta_{p^a})=\f1{\varphi(p^a)}.$$
By this and the induction hypothesis, for any $k=0,\ldots,n$ we have
$$\align&\ord_p\(C_{p^{a-1}}(k,r)\sum_{j=0}^{p-1}\zeta_{p^a}^{j(k-r)}
\l(1-\zeta_{p^a}^j\r)^{n-k}\)
\\\gs&\max\bg\{0,\,\l\lfloor\f{k-p^{a-2}}{\varphi(p^{a-1})}\r\rfloor\bg\}+\f{n-k}{\varphi(p^a)}
\\=&\max\bg\{0,\,\f{pk-p^{a-1}}{\varphi(p^a)}-\l\{\f{k-p^{a-2}}{\varphi(p^{a-1})}\r\}\bg\}
+\f{n-k}{\varphi(p^a)}
\\=&\max\bg\{\f{n-k}{\varphi(p^a)},\,\f{n-p^{a-1}}{\varphi(p^a)}+\f k{p^{a-1}}
-\l\{\f{k-p^{a-2}}{\varphi(p^{a-1})}\r\}\bg\}>\f{n-p^{a-1}}{\varphi(p^a)}.
\endalign$$
(Note that if $k\gs p^{a-1}$ then $k/p^{a-1}\gs1>\{(k-p^{a-2})/\varphi(p^{a-1})\}$.)
Therefore, from (5.2) we get that
$$\ord_p(C_{p^a}(n,r))>\f{n-p^{a-1}}{\varphi(p^a)}-1
\gs\l\lfloor\f{n-p^{a-1}}{\varphi(p^a)}\r\rfloor-1.$$
So $F_{p^a}(n,r)=(-p)^{-\lfloor(n-p^{a-1})/\varphi(p^a)\rfloor}C_{p^a}(n,r)\in\Z$
as desired. \qed

\medskip
\noindent{\it Proof of Theorem 1.4}.
(i) Write $n+d=p^{a-1}-1+m\varphi(p^a)$ with $m\in\N$.
Then, for any $k=0,\ldots,d$ we have
$$\l\lfloor\f{n+k-p^{a-1}}{\varphi(p^a)}\r\rfloor
=\l\lfloor m-\f{d-k+1}{\varphi(p^a)}\r\rfloor =m-1.$$

Below we use induction on $d$ to show the desired congruence (1.21).

In the case $d=0$ (i.e., $n-p^{a-1}\eq-1\ (\mo\ \varphi(p^a))$),
we have $F_{p^a}(n,r)\eq F_{p^a}(n,0)\ (\mo\ p)$
because
$$F_{p^a}(n,i)-F_{p^a}(n,i-1)=(-p)^{-m+1}C_{p^a}(n+1,i)
=-pF_{p^a}(n+1,i)$$
for all $i\in\Z$. Furthermore,
by a result of Weisman [We] (see also [SW, Theorem 1.5]),
$F_{p^a}(n,r)\eq1\ (\mo\ p)$ if $d=0$.

Now let $d>0$ and assume that the desired result holds for smaller values of $d$.
Clearly,  $(n+1)+(d-1)=p^{a-1}-1+m\varphi(p^a)$ and
$$\l\lfloor\f{n+1+k-p^{a-1}}{\varphi(p^a)}\r\rfloor=m-1\quad \t{for}\ \ k=0,\ldots,d-1.$$

If $r\gs0$ then
$$C_{p^a}(n,r)-C_{p^a}(n,0)=\sum_{0<i\ls r}(C_{p^a}(n,i)-C_{p^a}(n,i-1))
=\sum_{0<i\ls r}C_{p^a}(n+1,i);$$
if $r<0$ then
$$\align C_{p^a}(n,r)-C_{p^a}(n,0)=&\sum_{r<i\ls0}(C_{p^a}(n,i-1)-C_{p^a}(n,i))
\\=&-\sum_{r<i\ls0}C_{p^a}(n+1,i).
\endalign$$
Therefore
$$F_{p^a}(n,r)-F_{p^a}(n,0)=\cases\sum_{0<i\ls r}F_{p^a}(n+1,i)&\t{if}\ r\gs0,
\\-\sum_{r<i\ls 0}F_{p^a}(n+1,i)&\t{if}\ r<0.\endcases$$
By the induction hypothesis, whenever $i\in\Z$ we have
$$F_{p^a}(n+1,i)\eq\sum_{k=0}^{d-1}\bi{i+k-1}kF_{p^a}(n+1+k,0)\ (\mo\ p).$$
For any $k=0,\ldots,d-1$, if $r\gs0$
$$\sum_{0<i\ls r}\bi{i+k-1}k=\sum_{j=0}^{r+k-1}\bi jk=\bi{r+k}{k+1}$$
by an identity of S.-C. Chu (cf. [GKP, (5.10)]); if $r<0$ then
$$\align -\sum_{r<i\ls0}\bi{i+k-1}k=&(-1)^{k+1}\sum_{r<i\ls0}\bi{-i}k
=(-1)^{k+1}\sum_{j=0}^{-r-1}\bi jk
\\=&(-1)^{k+1}\bi{-r}{k+1}=\bi{r+k}{k+1}.
\endalign$$
Thus, by the above, $F_{p^a}(n,r)$ is congruent to
$$F_{p^a}(n,0)+\sum_{k=0}^{d-1}\bi{r+k}{k+1}F_{p^a}(n+1+k,0)
=\sum_{k=0}^d\bi{r+k-1}kF_{p^a}(n+k,0)$$
modulo $p$. This concludes the induction proof of (1.21). \qed

(ii) In the case $a=1$,
the desired results in Theorem 1.4(ii) follow from Corollaries 1.3 and 1.7.

Now we let $a\gs2$ and $r\in\Z$.
Write $n=p^{a-2}(pn_1+n_0)+s$ and $r=p^{a-2}(pr_1+r_0)+t$,
where $s,t\in\{0,\ldots,p^{a-2}-1\}$, $n_0,r_0\in\{0,\ldots,p-1\}$ and $n_1\in\N$ and $r_1\in\Z$.

If $p^{a-1}\ls n<p^a$, then
$$F_{p^a}(n,r)=C_{p^a}(n,r)=\bi n{\{r\}_{p^a}}(-1)^{\{r\}_{p^a}},$$
and in particular $\ord_p(C_{p^a}(n,0))=0=\lfloor(n-p^{a-1})/\varphi(p^a)\rfloor$.

Below we assume that $n\gs 2p^{a-1}$ (i.e., $n_1\gs2$). By [SD, Theorem 1.7],
$$F_{p^a}(n,r)\eq(-1)^t\bi st F_{p^2}(pn_1+n_0,pr_1+r_0)\ (\mo\ p).$$

If $p\mid n_1$, or $p-1\nmid n_1-1$, or $n_0=r_0=p-1$,
then by [SW, Theorem 1.2] in the case $l=0$, we have
$$F_{p^2}(pn_1+n_0,pr_1+r_0)\eq(-1)^{r_0}\bi{n_0}{r_0}F_p(n_1,r_1)\ (\mo\ p)$$
and hence $F_{p^a}(n,r)\eq b_{n,r}F_p(n_1,r_1)\ (\mo\ p)$, where
$$\align b_{n,r}:=&(-1)^{\{r\}_{p^{a-1}}}\bi{\{n\}_{p^{a-1}}}{\{r\}_{p^{a-1}}}
=(-1)^{p^{a-2}r_0+t}\bi{p^{a-2}n_0+s}{p^{a-2}r_0+t}
\\\eq&(-1)^{t}\bi st(-1)^{r_0}\bi{n_0}{r_0}\ (\mo\ p) \ \  (\t{by Lucas' theorem (cf. [HS])}).
\endalign$$
By Corollary 1.3, there is an $r_1'\in\Z$ such that $F_p(n_1,r_1')\not\eq0\ (\mo\ p)$.
Thus, if $p\mid n_1$ or $p-1\nmid n_1-1$, then
$$F_{p^a}(n,p^{a-1}r_1')\eq F_p(n_1,r_1')\not\eq0\ (\mo\ p).$$
If $n_0=p-1$, then
$$F_{p^a}(n,p^{a-2}(pr_1'+p-1))
\eq (-1)^{p-1}\bi{p-1}{p-1}F_p(n_1,r_1')\not\eq0\ (\mo\ p).$$

When $p\nmid n_1$, $p-1\mid n_1-1$ and $n_0<r_0$, by applying the second part of [SW, Theorem 1.2]
in the case $l=0$, we have
$$F_{p^2}(pn_1+n_0,pr_1+r_0)\eq[\![n_1>1]\!]\f{(-1)^{n_0}n_1}{r_0\bi{r_0-1}{n_0}}
=\f{(-1)^{n_0}n_1}{r_0\bi{r_0-1}{n_0}}\ (\mo\ p)$$
and hence
$$F_{p^a}(n,r)\eq(-1)^{n_0+t}\f{n_1\bi st}{r_0\bi{r_0-1}{n_0}}\ (\mo\ p).$$
In particular, if $p\nmid n_1$, $p-1\mid n_1-1$ and $n_0<p-1$, then
$$F_{p^a}(n,p^{a-2}(n_0+1))\eq\f{(-1)^{n_0}n_1}{n_0+1}\not\eq0\ (\mo\ p).$$

In view of the above, we already have (1.22).

To prove the congruence in (1.23), we should also
consider the case $p\nmid n_1$, $p-1\mid n_1-1$ and $n_0\gs r_0$.
By [SW, Lemmas 3.2 and 3.3],
$$\align &p^{-\lfloor(pn_1+n_0-p)/\varphi(p^2)\rfloor}C_{p^2}(pn_1+n_0,pr_1+r_0)
\\&-(-1)^{r_0}\bi{n_0}{r_0}p^{-\lfloor(n_1-1)/(p-1)\rfloor}C_p(n_1,r_1)
\\\eq&(-1)^{n_1-1}p^{-\lfloor(n_1-1-1)/(p-1)\rfloor}C_p(n_1-1,r_1)
(-1)^{n_1+r_0}n_1\bi{n_0}{r_0}\f{\sigma_{n_0,r_0}(n_1)}p
\\\eq&-(-1)^{r_0}\bi{n_0}{r_0}p^{-(n_1-1)/(p-1)+1}C_p(n_1-1,r_1)
n_1\f{\sigma_{n_0,r_0}(n_1)}p\ (\mo\ p),
\endalign$$
where
$$\sigma_{n_0,r_0}(n_1)=1+(-1)^p\f{\prod_{1\ls i\ls p,\,i\not=p-r_0}(p(n_1-1)+r_0+i)}
{\prod_{1\ls i\ls p,\, i\not=p-(n_0-r_0)}(n_0-r_0+i)}\eq0\ (\mo\ p).$$
Therefore
$$\align &F_{p^2}(pn_1+n_0,pr_1+r_0)-(-1)^{r_0}\bi{n_0}{r_0}F_p(n_1,r_1)
\\\eq&(-1)^{r_0}\bi{n_0}{r_0}F_p(n_1-1,r_1)n_1\f{\sigma_{n_0,r_0}(n_1)}p\ (\mo\ p)
\endalign$$
and hence
$$F_{p^a}(n,r)\eq b_{n,r}\(F_p(n_1,r_1)+F_p(n_1-1,r_1)n_1\f{\sigma_{n_0,r_0}(n_1)}p\)\ (\mo\ p),$$

Observe that $n+p^a(p-1)=p^{a-2}(pn_1'+n_0)+s$ with $n'_1=n_1+p(p-1)$.
Clearly $F_p(n_1',r_1)\eq F_p(n_1,r_1)\ (\mo\ p)$ by Corollary 1.7,
and $\sigma_{n_0,r_0}(n_1')\eq\sigma_{n_0,r_0}(n_1)\ (\mo\ p^2)$ if $n_0\gs r_0$.
Thus, by the above, $F_{p^a}(n+p^a(p-1),r)\eq F_{p^a}(n,r)
\ (\mo\ p)$. This concludes the proof. \qed

\enddocument